\documentclass[12pt]{article}
\usepackage[reqno]{amsmath} 
\usepackage{amssymb, amsthm}
\usepackage{enumerate}  
\usepackage{url} 
\usepackage{hyperref}
\usepackage{graphicx}

\newtheoremstyle{plainsl}%
	{\topsep}
	{\topsep}
	{\slshape} 
	{}
	{\normalfont\bfseries}
	{.}
	{ }
	{}

\swapnumbers

\theoremstyle{plainsl}
\newtheorem{theorem}{Theorem}[section]
\newtheorem{lemma}[theorem]{Lemma}
\newtheorem{corollary}[theorem]{Corollary}

\renewcommand\proof{\noindent\textsl{Proof. }}
\newcommand\sqr[2]{{\vbox{\hrule height.#2pt
    \hbox{\vrule width.#2pt height#1pt \kern#1pt
        \vrule width.#2pt}\hrule height.#2pt}}}
\renewcommand\qed{%
	\ifmmode\eqno\sqr53
	\else\nolinebreak\ \hfill\sqr53\medbreak\fi}

\numberwithin{equation}{section}

\newcommand{\ii}{\mathrm{i}}
\newcommand{\dd}{d}

\newcommand{\extremal}{spectrally extremal\ }

\newcommand{\e}{\mathrm{e}}

\title{Spectrally extremal vertices, strong cospectrality and state transfer}
\author{Gabriel Coutinho\footnote{Dep. of Combinatorics and Optimization, University of Waterloo. Research partialy funded by the CAPES foundation, Ministry of Education, Brazil. \texttt{gcoutinho@uwaterloo.ca}}}

\begin{document}

\maketitle

\begin{abstract}
In order to obtain perfect state transfer between two sites in a network of interacting qubits, their corresponding vertices in the underlying graph must satisfy a combinatorial property called strong cospectrality. Here we determine the structure of graphs containing pairs of vertices which are strongly cospectral and satisfy a certain extremal property related to the spectrum of the graph. If the graph satisfies this property globally and is regular, we also show that the existence of a partition of the vertex set into pairs of vertices at maximum distance admitting perfect state transfer forces the graph to be distance-regular.
\end{abstract}

\section{Introduction}

Let $X$ be a simple undirected graph, and consider its $01$-adjacency matrix $A = A(X)$. For $u \in V(X)$, we denote by $e_u$ the vector of the canonical basis corresponding to $u$ in the ordering of the rows of $A$. The matrix operator
\[U(t) = \exp(\ii t A),\]
defined for every real $t \geq 0$, represents a \textsl{continuous-time quantum walk} on $X$. It was shown in \cite{ChristandlPSTQuantumSpinNet} that the dynamics of quantum state transfer in a network of interacting qubits in the $XY$-coupling model is determined by properties of $U(t)$. More specifically, we will say that $X$ admits \textsl{perfect state transfer} between vertices $u$ and $v$ at time $\tau$ if
\[|e_v^T U(\tau) e_u| = 1.\]
The problem of determining which graphs admit perfect state transfer has been studied recently in a good number of papers. For example, it was solved for paths and hypercubes (see \cite{ChristandlPSTQuantumSpinNet2}), circulant graphs (see \cite{BasicPetkovicPSTCirculant}) and cubelike graphs (see \cite{GodsilCheungPSTCubelike}). The effect of certain graph operations was considered in \cite{Angeles-CanulPSTcirculant} and in \cite{TamonPSTQuotient}, and some recent surveys are found in \cite{KendonTamon} and \cite{GodsilStateTransfer12}.

The study of state transfer or more generally the study of continuous-time quantum walks on graphs can be divided into two aspects. The combinatorial properties that must be satisfied by the graph, and the number theoretic properties that must be satisfied by the eigenvalues of the graph. In the case of perfect state transfer between vertices $u$ and $v$, the key combinatorial property is that $u$ and $v$ must be \textsl{strongly cospectral}, that is, the projections of $e_u$ and $e_v$ in each eigenspace of the graph must be parallel with the same magnitude.

In this paper, we will examine this property in the context of graphs which are extremal with respect to the lower bound on the number of eigenvalues given by the diameter plus one. State transfer on such graphs was considered in \cite{ZhouBuShen}. Our approach is nevertheless more general and results in deeper consequences.

We will consider the extremal spectral property locally, and show that a pair of strongly cospectral extremal vertices at maximum distance must be singletons in a pseudo equitable distance partition. This will be our key intermediate step to show that if an extremal regular graph of diameter $d$ can be partitioned into pairs of vertices at distance $d$ such that perfect state transfer happens in each pair, then the graph is distance-regular. We will also characterize precisely which number theoretic conditions must be satisfied by the eigenvalues of the graph in order to achieve perfect state transfer between strongly cospectral vertices. In the context of distance-regular graphs, this will provide an alternate elementary proof of \cite[Corollary 4.5]{CoutinhoGodsilGuoVanhove}. 

\section{Preliminaries}

Let $X$ be a graph, $u,v \in V(X)$, and $A$ the adjacency matrix of $X$. We denote the spectral decomposition of $A$ into orthogonal projections by
\[A = \sum_{r = 0}^d \theta_r E_r.\]

Vertices $u$ and $v$ are called \textsl{cospectral} if the adjacency matrices of $X\backslash u$ and $X\backslash v$ have the same spectrum. 
\begin{lemma}[\cite{GodsilAlgebraicCombinatorics}, Chapter 4, Lemma 1.1] \label{lem:0.1}
The following are equivalent.
\begin{enumerate}[(i)]
\item Vertices $u$ and $v$ are cospectral.
\item For all $k \in \mathbb{Z}^+$,
\[(A^k)_{u,u} = (A^k)_{v,v}.\]
\item For all $r=0,...,d$,
\[(E_r)_{u,u} = (E_r)_{v,v}.\]
\end{enumerate} \ \qed
\end{lemma}

Recall that we denote the characteristic vector of a vertex $w$ by $e_w$. Vertices $u$ and $v$ are \textsl{strongly cospectral} if and only if $E_r e_u = \pm E_r e_v$ for all $r = 0,...,d$. Thus strongly cospectral vertices must be cospectral.

We define the \textsl{eigenvalue support} of $u$ as the set of all eigenvalues $\theta_r$ such that $E_r e_u \neq 0$, and we denote it by $\Phi_u$. We say that the \textsl{dual degree} of $u$ is $d^*(u) = |\Phi_u| - 1$. The \textsl{eccentricity} (or covering radius) of $u$, denoted by $\varepsilon_u$, is the maximum distance between $u$ and any other vertex of the graph.

The \textsl{distance partition} relative to $u \in V(X)$ is the partition of the vertex set of $X$ in which each class is formed by vertices at a fixed distance from $u$. A partition of the vertex set of a graph is called \textsl{equitable} if the number of edges between a vertex $v$ and a class $C$ of the partition depends only on the class containing $v$. The collection of all these numbers are the \textsl{parameters} of the equitable partition.

If $D$ is a non-singular diagonal matrix, we consider the matrix $(D^{-1} A D)$. This can be seen as the adjacency matrix of a weighted directed graph, where $(D^{-1}AD)_{uv}$ indicates the weight of the directed arc from $u$ to $v$. We will call a partition of $V(X)$ \textsl{equitable with respect to $D$} if the sum of the weights in $(D^{-1} A D)$ from a vertex $v$ to a class $C$ of the partition depends only of the class containing $v$.

When $D$ is chosen to be the diagonal matrix whose entries are taken from the eigenvector of $A$ corresponding to the largest eigenvalue, an equitable partition with respect to $D$ will simply be called \textsl{pseudo equitable}. Note that for regular graphs, a partition is pseudo equitable if and only if it is equitable in the original sense.

A graph is called \textsl{distance-regular} if the distance partition relative to each vertex is equitable and the parameters do not depend on the vertex. A distance-regular graph of diameter $d$ is called \textsl{antipodal} if the relation defined on the vertex set by having two vertices related if they are at distance $0$ or $d$ is an equivalence relation.

\begin{theorem}[\cite{GodsilShawe-Taylor}, Theorem 2.2] \label{thm:0}
If $X$ is a regular graph such that the distance partition of every vertex is equitable, then $X$ is distance-regular.
\end{theorem}

The \textsl{walk module} of a vertex $u$ in a graph $X$ is the subspace
\[W_u = \langle \{ A^i e_u\}_{i \geq 0} \rangle.\]
It follows that
\[W_u = \langle \{ E_r e_u \}_{\theta_r \in \Phi_u} \rangle,\]
and because the vectors $\{A^i e_u\}_{i = 0}^{\varepsilon_u}$ are all independent, we have that
\begin{align}\label{eq:1}\varepsilon_u \leq d^*(u).\end{align}

If equality holds, we will say that $u$ is a \textsl{\extremal} vertex. The work of Fiol, Garriga and others provides a vast literature about graphs with \extremal vertices and associated concepts (see for instance \cite{FiolGarrigaSpectralExcessTheorem}, \cite{FiolGarrigaDalfoVanDamAlmostDRGs}, \cite{LeeWengSpectralExcess} and \cite{vanDamSpectralExcess}). Here we will use the following result.

\begin{theorem}[\cite{FiolGarrigaYebra}, Theorem 6.3] \label{thm:1}
The distance partition relative to a vertex $u$ is pseudo equitable if and only if $u$ is \extremal and there exists a polynomial $p(x)$ such that $p(A)e_u$ is a $01$-vector whose support are the vertices at distance $\varepsilon_u$ from $u$.
\end{theorem}

On the topic of perfect state transfer, we will use the following two results.

\begin{theorem}[\cite{GodsilPerfectStateTransfer12}, Theorem 6.1]\label{thm:2}
If $X$ admits perfect state transfer between $u$ and $v$, then either all eigenvalues in the eigenvalue support of $u$ are integers, or they are all of the form $\theta_r = \frac{1}{2} (a + b_r \sqrt{\Delta})$ where $\Delta$ is a square-free integer and $a$ and $b_r$ are integers for all $r$.
\end{theorem}

\begin{lemma}[\cite{GodsilStateTransfer12}, Lemma 11.1] \label{lem:2}
If $X$ admits perfect state transfer between $u$ and $v$, then $u$ and $v$ are strongly cospectral.
\end{lemma}

\section{Strong cospectrality on \extremal vertices}

It turns out that when two strongly cospectral vertices are also spectrally extremal, each of them is a singleton in the distance partition of the other.

\begin{lemma} \label{lem:1}
Let $u,v \in V(X)$, with $g = d(u,v)$. Suppose $u$ is a \extremal vertex. If $u$ and $v$ are strongly cospectral, then the following hold.
\begin{enumerate}[(i)]
\item If $d(u,w) = d(u,v)$, then $w=v$.
\item If $z \in V(X)$, $\Phi_z = \Phi_u$ and if $d(z,w) = d(u,v)$ for some $w \in V(X)$, then $(A^g)_{z,w} \leq (A^g)_{u,v}$. Equality occurs if and only if $z$ and $w$ are strongly cospectral.
\end{enumerate} 
\end{lemma}
\proof
Suppose $\Phi_u = \{\theta_0,\theta_1,....,\theta_{d^*}\}$. For all $r \in \{0,...,d^*\}$, let $\sigma_r \in \{+1,-1\}$ be such that
\[E_r e_v = \sigma_r E_r e_u.\]
Let $p(x)$ be the polynomial of minimum degree satisfying $p(\theta_r) = \sigma_r$ for all $r$. Then it follows that
\[p(A) e_u = e_v.\]
Because $\varepsilon_u = d^*$, $p(A) e_u$ must be non-zero on the entries corresponding to vertices whose distance to $u$ is the degree of $p(x)$. Hence $\textrm{deg}\  p(x) = g$, and $v$ is the unique vertex at distance $g$ from $u$.

To see (ii), first note that $\langle p(A) e_z ,  p(A) e_z \rangle = 1$, so the absolute value of each entry in $p(A)e_z$ is at most $1$. Let $p(x) = a_gx^g + ... + a_0$. Then $p(A)e_u = e_v$ implies that
\[a_g = \frac{1}{(A^g)_{u,v}},\]
and thus
\[1 \geq |p(A)_{z,w}| = a_g (A^g)_{z,w} = \frac{(A^g)_{z,w}}{(A^g)_{u,v}}.\]
\qed

We are also going to need the following characterization of strongly cospectral vertices.

\begin{lemma}\label{lem:5}
Let $u,v \in V(X)$. The following are equivalent.
\begin{enumerate}[(i)]
\item Vertices $u$ and $v$ are cospectral, and there exists a polynomial $p(x)$ such that ${p(A) e_u = e_v}$.
\item The vertices $u$ and $v$ are strongly cospectral.
\end{enumerate}
Moreover, if $u$ and $v$ are cospectral, then any polynomial satisfying $p(A) e_u = e_v$ is such that $p(A)e_v = e_u$ and $p(\theta_r) = \pm 1$ for all $\theta_r \in \Phi_u$.
\end{lemma}
\proof
The implication (ii) $\implies$ (i) is trivial. To see the converse, let $p(x)$ be a polynomial satisfying $p(A) e_u = e_v$. Because $p(A)$ is a symmetric matrix, it follows that $(p(A)^{2})_{u,u} = 1$. Vertices $u$ and $v$ are cospectral, so Theorem \ref{lem:0.1} implies that $(p(A)^2)_{v,v} = 1$. Thus $p(A)e_v$ is a unitary vector, but $p(A)_{u,v} = 1$, implying that $p(A)e_v = e_u$. As a consequence, $p(A)^2 e_u = e_u$, and so if $\theta_r \in \Phi_u$, it follows that $p(\theta_r) = \pm 1$. This shows that $u$ and $v$ are strongly cospectral.
\qed

Here we introduce a definition. We say that $u$ and $v$ are (a pair of) \textsl{antipodal\footnote{The name \textit{antipodal} has been used in different contexts. We are consistent with at least one its uses, namely, if a regular graph can be partitioned into pairs of antipodal vertices, then it is an antipodal distance regular graph.} vertices} if 
\begin{itemize}
\item there is a pseudo equitable partition which is simultaneously the distance partition of $u$ and $v$,
\item $\{u\}$ and $\{v\}$ are singletons in this partition at maximum distance from each other, and
\item $u$ and $v$ are cospectral.
\end{itemize}

If the graph is regular, then a pair $u$ and $v$ of antipodal vertices is cospectral if and only if the parameters of the distance partition of $u$ are equal to the parameters of the distance partition of $v$. This is a consequence of Lemma \ref{lem:0.1} and of the fact that the number of closed walks of any length on vertices whose distance partition is equitable is determined by the parameters of the partition. 

\begin{theorem} \label{thm:4}
If $u$ and $v$ are antipodal vertices in $X$, then $u$ and $v$ are \extremal vertices and they are strongly cospectral. On the other hand, if $u$ is \extremal, $u$ and $v$ are strongly cospectral, and their distance is equal to their eccentricity, then $u$ and $v$ are antipodal vertices.
\end{theorem}
\proof
If $u$ and $v$ are antipodal, then the weaker direction of Theorem \ref{thm:1} implies that $u$ is \extremal and that there is a polynomial $p(x)$ such that
\[p(A) e_u = e_v.\]
From Lemma \ref{lem:5}, we have that $u$ and $v$ are strongly cospectral.

The converse is an immediate consequence of Theorem \ref{thm:1} and the fact that strongly cospectral vertices are cospectral.
\qed

We would like to drop the condition on the theorem above that requires $u$ and $v$ to be at maximal distance. In other words, we would like to believe that a pair of spectrally extremal strongly cospectral vertices must be at maximal distance from each other. For instance, this is true for $2$-connected graphs. In particular, Lemma \ref{lem:1} implies that $v$ is a cut-vertex of $X$ unless $v$ is at maximal distance from $u$. So if $X$ is $2$-connected, it follows that $u$ and $v$ must be at maximal distance. 

The following lemma is a step towards extending this last observation to all graphs. We will use it to derive some interesting consequences.

\begin{lemma} \label{lem:3}
Suppose $u$ is a \extremal vertex of $X$, and suppose $u$ and $v$ are strongly cospectral. Let $p(x)$ be the polynomial satisfying $p(A)e_u = e_v$, with $p(\theta_r) = \sigma_r \in \{+1,-1\}$ for all $\theta_r \in \Phi_u$. Let $X'$ be the component of $X \backslash v$ containing $u$. Then $p(x)$ is the minimal polynomial with respect to $u$ in $X'$ (up to a constant).
\end{lemma}
\proof
Let $d(u,v) = g$ and $A'=A(X')$. From Lemma \ref{lem:1} (i), we have that $v$ is the unique vertex at distance $g$ from $u$. Note that walks of length $g$ pass by $v$ only if $v$ is its final vertex, so the entries of $p(A)e_u$ relative to vertices at distance at most $g-1$ from $u$ are equal to the respective entries of $p(A')e_u$, thus $p(A')e_u = 0$. Because the eccentricity of $u$ in $X'$ is $g-1$ and $p(x)$ has degree $g$, it follows that it is the minimal polynomial up to a constant.
\qed

\begin{corollary}
Let $u,v \in V(X)$. Suppose $\Phi_u = \{\theta_0,....,\theta_{\dd^*}\}$, ordered in such way that $\theta_r > \theta_{r+1}$ for all $r$. If $u$ and $v$ are spectrally extremal and strongly cospectral, and $p(x)$ is such that $p(A) e_u = e_v$, then there is no index $r \in \{0,...,\dd^*\}$ such that
\[p(\theta_r) = p(\theta_{r+1}) = p(\theta_{r+2}).\]
\end{corollary}
\proof
Suppose otherwise that there is such index, say $s$. From Lemma \ref{lem:3}, the roots of $p(x)$ are the eigenvalues of $X \backslash v$ in the support of $u$. Interlacing (see \cite[Theorem 2.5.1]{BrouwerHaemers}) implies that there are no two roots of $p(x)$ between $\theta_r$ and $\theta_{r+1}$ for any $r$, hence $p(\theta_s) = p(\theta_{s+1}) = p(\theta_{s+2})$ implies that there are three roots of $\frac{\dd}{\dd x}p(x)$ between two of its real roots. This is a contradiction to the fact that all roots of $p(x)$ are real.
\qed

If we know that the pair of strongly cospectral vertices is at maximal distance, we can actually determine the values of $p(\theta_r)$ for all $r$.

\begin{theorem} \label{lem:4}
Let $u,v \in V(X)$. Suppose $\Phi_u = \{\theta_0,....,\theta_{\dd^*}\}$, ordered in such way that $\theta_r > \theta_{r+1}$ for all $r$. Then $u$ and $v$ are antipodal if and only if, for all $r \in \{0,...,\dd^*\}$,
\[E_r e_v = (-1)^r E_r e_u.\]
\end{theorem}
\proof
Let $d = d(u,v)$, and $p(x)$ the polynomial of degree $d$ such that $p(A) e_u = e_v$. If $p(x)$ is such that $p(\theta_r) = (-1)^r$, then $p(x)$ has at least $\dd^*$ roots, so $d \geq \dd^*$, and hence it could only be that $d = \dd^*$. So $u$ and $v$ are spectrally extremal, strongly cospectral, and at maximal distance. It follows from Theorem \ref{thm:4} that they are a pair of antipodal vertices.

For the converse, let $p(x)$ be the polynomial satisfying $p(\theta_r) = \sigma_r \in \{+1,-1\}$ with $p(A)e_u = e_v$. The existence of such $p(x)$ is implied in Theorem \ref{thm:4}. Let $q(x)$ be the polynomial of minimal degree that satisfies $q(\theta_r) = (-1)^r$ for $r \in \{0,...,\dd^*\}$. Therefore
\begin{align*}
1 & \geq |q(A)_{u,v}| \\
& = \left| \sum_{r=0}^{\dd^*} (-1)^r \prod_{s \neq r} \frac{1}{\theta_r - \theta_s}\right|(A^{\dd^*})_{u,v} \\
& = \left( \sum_{r=0}^{\dd^*} (-1)^r \prod_{s \neq r} \frac{1}{\theta_r - \theta_s}\right)   (A^{\dd^*})_{u,v}, \quad \text{because all terms are positive,}\\
& \geq \left( \sum_{r=0}^{\dd^*} \sigma_r \prod_{s \neq r} \frac{1}{\theta_r - \theta_s} \right)  (A^{\dd^*})_{u,v} \\
& = p(A)_{u,v}\\
& = 1.
\end{align*}
Note that equality holds throughout if and only if $\sigma_r = (-1)^r$, as we wanted.
\qed

If a graph $X$ has diameter $d$, then Equation \ref{eq:1} implies that $X$ has at least $d+1$ distinct eigenvalues. We say that $X$ is \textit{\extremal} if equality holds. Note that every \extremal graph contains at least one pair of \extremal vertices.

\begin{theorem} \label{thm:3}
Suppose $X$ is a \extremal regular graph on $n$ vertices of diameter $d$, and that its distinct eigenvalues are $\theta_0>...>\theta_d$. Suppose $u$ and $v$ are vertices at distance $d$. Then $u$ and $v$ are antipodal if and only if
\[n \prod_{s = 0}^d \frac{1}{\theta_0 - \theta_s } = \sum_{r = 0}^d (-1)^r \prod_{s \neq r} \frac{1}{\theta_r - \theta_s}.\]
\end{theorem}
\proof
Let $p(x)$ be a polynomial such that $p(A) = E_0$, and let $J$ denote the all $1$s matrix. Because the graph is regular $E_0 = (1/n) J$, and so if $p(x) = a_dx^d + ... + a_0$, it follows that
\[(A^d)_{u,v} = \frac{1}{n} \left( \prod_{s=0}^d \frac{1}{\theta_0 - \theta_s}\right)^{-1} \]
for all vertices $u$ and $v$ at distance $d$. The result now follows from Lemma \ref{lem:4}.
\qed

\begin{corollary}\label{cor:3.8}
Suppose $X$ is a \extremal regular graph on $n$ vertices of diameter $d$. If the eccentricity of every vertex is $d$ and if
\[n \prod_{s = 0}^d \frac{1}{\theta_0 - \theta_s } = \sum_{r = 0}^d (-1)^r \prod_{s \neq r} \frac{1}{\theta_r - \theta_s},\]
then $X$ is an antipodal distance regular graph.
\end{corollary}
\proof
It follows from Theorem \ref{thm:0} and Theorem \ref{thm:3}.
\qed

\section{State transfer on \extremal graphs}

In this section, we relate the work in the past section to state transfer.

\begin{theorem}\label{thm:5}
Suppose $u$ is a \extremal vertex of $X$, and $v$ is a vertex at maximal distance from $u$. Let $\Phi_u = \{\theta_0,...,\theta_{\dd^*}\}$, where, for some integers $a$, $\Delta$ and $b_r$ with $r = 0,...,\dd^*$, we have $\theta_r = \frac{a + b_r\sqrt{\Delta}}{2}$ satisfying $\theta_r > \theta_{r+1}$. Then $X$ admits perfect state transfer between vertices $u$ and $v$ if and if only if the following conditions hold.
\begin{enumerate}[(i)]
\item The vertices $u$ and $v$ are antipodal.
\item If $r$ is odd, then the powers of two appearing in the factorization of each of the differences $b_0 - b_r$ are constant, let us say $\alpha$.
\item If $r$ is even, then the power of two in the factorization of each $b_0 - b_r$ is larger than $\alpha$.
\end{enumerate}
If the conditions hold, perfect state transfer happens at time $\frac{\pi}{2^\alpha}$ (and any of its odd multiples).
\end{theorem}
\proof
From Theorem \ref{thm:2}, perfect state transfer implies strong cospectrality. The vertices are at maximal distance from each other, and if they are strongly cospectral, Theorem \ref{thm:4} says that they are antipodal. Note that
\[U(t) e_u = \sum_{r = 0}^{\dd^*} \e^{\ii t \theta_r} E_r e_u.\]
Lemma \ref{lem:4} implies that
\[e_u = \sum_{r=0}^{ \dd^* } (-1)^r E_r e_v,\]
therefore perfect state transfer is now equivalent to existing a time $t$ such that
\[\frac{\e^{\ii t \theta_0}}{\e^{\ii t \theta_r}} = (-1)^r,\]
and this is equivalent to
\[t (\theta_0 - \theta_r) = k_r \pi,\]
where $k_r$ is an integer with the same parity as $r$. This condition is equivalent to (ii) and (iii), and also gives the expression for the time.
\qed

\begin{corollary}
Suppose $X$ is a \extremal graph of diameter $d$ and eigenvalues $\theta_0 > ... > \theta_d$. Suppose $X$ admits perfect state transfer at time $\tau$ between vertices $u$ and $v$ at distance $d$. If $d(z,w) = d$, then $(A^d)_{z,w} \leq (A^d)_{u,v}$, and equality happens if and only if $X$ admits perfect state transfer between $z$ and $w$ at time $\tau$.
\end{corollary}
\proof
The inequality is a straightforward application of Lemma \ref{lem:1} (ii). If equality happens, then Lemma \ref{lem:1} (ii) says that $z$ and $w$ are strongly cospectral. They are at maximum distance, so by Theorem \ref{thm:4}, they are antipodal. Combining this with Theorem \ref{thm:5} and the fact that the eigenvalue support of $z$ and $w$ is equal to the eigenvalue support of $u$ and $v$, it follows that perfect state transfer between $z$ and $w$ happens at time $\tau$.
\qed

In the case where $X$ is a regular graph, we can say more.

\begin{corollary}
Suppose $X$ is a \extremal regular graph on $n$ vertices of diameter $d$, distinct eigenvalues $\theta_0>...>\theta_d$. Then $X$ admits perfect state transfer between any pair $(u,v)$ of vertices at distance $d$ if and only if
\begin{enumerate}[(i)]
\item All eigenvalues are integers.
\item If $r$ is odd, then the powers of two appearing in the factorization of each of the differences $\theta_0 - \theta_r$ are constant, let us say $\alpha$.
\item If $r$ is even, then the power of two in the factorization of $\theta_0 - \theta_r$ is larger than $\alpha$.
\item The following equality holds
\[n \prod_{s = 0}^d \frac{1}{\theta_0 - \theta_s } = \sum_{r = 0}^d (-1)^r \prod_{s \neq r} \frac{1}{\theta_r - \theta_s}.\]
\end{enumerate}
\end{corollary}
\proof
Because the graph is spectrally extremal, all eigenvalues are in the eigenvalue support of the vertices at distance $d$. Note that the largest eigenvalue is integer, so Theorem \ref{thm:2} implies that all eigenvalues are integers. Thus it follows from Theorems \ref{thm:3} and \ref{thm:5} that the other conditions are equivalent to perfect state transfer between vertices at distance $d$.
\qed

The following corollary is an immediate consequence of the result above and Corollary \ref{cor:3.8}.

\begin{corollary}
If $X$ is a \extremal regular graph of diameter $d$ in which the eccentricity of every vertex is $d$, and if perfect state transfer happens between any pair of vertices at distance $d$, then $X$ is a distance-regular graph.
\end{corollary}

We finally observe that the result above is a good example on how the existence of perfect state transfer can imply deep structural properties in a graph.

\section*{Acknowledgement}

I'd like to acknowledge Chris Godsil for very fruitful discussions about perfect state transfer and more specifically about the results in this paper. I also acknowledge M. A. Fiol, who was visiting Waterloo during the end of my PhD and with whom I discussed the subject of spectrally extremal graphs.

\bibliographystyle{plain}
\bibliography{qwalks}

\end{document}